\documentclass[conference]{IEEEtran}
\usepackage{multirow}
\usepackage{multicol}
\usepackage{lipsum}
\usepackage{graphicx}
\usepackage{graphics}
\usepackage{amsmath}
\usepackage{amsbsy}
\usepackage{blindtext}
\usepackage{tcolorbox}
\usepackage{amssymb}
\usepackage{scalerel}
\usepackage{mathabx}
\usepackage{verbatim}
\usepackage{booktabs}
\usepackage{rotating,tabularx}
\usepackage{mathrsfs} 
\usepackage{hyperref}

\usepackage{adjustbox}
\usepackage{soul}
\usepackage{xcolor}
\usepackage{cite}
\usepackage{tikz}
\usepackage{color, colortbl}
\usepackage[first=0,last=9]{lcg}
\definecolor{LightGray}{rgb}{0.7,0.7,0.7}

\usepackage[wby]{callouts}
\usetikzlibrary{shapes.multipart}

\usepackage{array}
\usepackage{makecell}

\allowdisplaybreaks

\usepackage{amsthm}
\theoremstyle{definition}

\theoremstyle{remark}

\usepackage[utf8]{inputenc}
\usepackage[english]{babel}

\usepackage[caption=false,font=footnotesize]{subfig}
\usepackage[T1]{fontenc}
\usepackage{scalerel,stackengine}
\newcommand\reallywidecheck[1]{%
\savestack{\tmpbox}{\stretchto{%
  \scaleto{%
    \scalerel*[\widthof{\ensuremath{#1}}]{\kern-.6pt\bigwedge\kern-.6pt}%
    {\rule[-\textheight/2]{1ex}{\textheight}}
  }{\textheight}%
}{0.5ex}}%
\stackon[1pt]{#1}{\scalebox{-1}{\tmpbox}}%
}
\stackMath
\IEEEoverridecommandlockouts
\IEEEoverridecommandlockouts

\definecolor{shadecolor}{RGB}{190,190,190}

\renewcommand{\baselinestretch}{0.967}

\newcommand*{\mrn}{\textcolor{black}}

\newif\ifarxiv
\arxivtrue
\arxivfalse

\begin{document}

\title{\LARGE\bf
Analyzing the Impact of AC False Data Injection Attacks on Power System Operation}


\author{\IEEEauthorblockN{1\textsuperscript{st} Mohammadreza Iranpour}
\IEEEauthorblockA{\textit{Department of Electrical and Computer Engineering} \\
\textit{California State University Northridge (CSUN)}\\
Los Angels, USA \\
\and
\IEEEauthorblockN{2\textsuperscript{nd} Mohammad Rasoul Narimani}
\IEEEauthorblockA{\textit{Department of Electrical and Computer Engineering} \\
\textit{California State University Northridge (CSUN)}\\
Los Angels, USA \\
Rasoul.narimani@csun.edu}}
\and

\thanks{This work has been supported from NSF contract \#2308498.}%
}

\maketitle

\begin{abstract}
\mrn{False Data Injection (FDI) attacks are a significant threat to modern power systems. Although numerous research studies have focused on FDI attacks on power systems, these studies have primarily concentrated on designing or detecting DC FDI attacks, with less attention given to the impact analysis of AC FDI attacks. AC FDI attacks are potentially more harmful as they can easily bypass bad data detection (BDD) algorithms. In this paper, we present a unified approach to investigate the impact of AC FDI attacks on power transmission lines using the PowerWorld simulator. We also investigate the impact of different FDI attack designs, including those optimally designed to evade BDD algorithms and compare them accordingly. Our findings demonstrate that in designing optimal AC FDI attacks, a trade-off between the residuals of state variables and the corresponding impacts of the proposed attack should be considered. This is because optimal attacks result in fewer changes in the attacked variable states and their estimated residuals compared to arbitrary AC FDI attacks. Moreover, the impacts of optimal AC FDI attacks can be less severe than those of arbitrary attacks. We implement and analyze the proposed approach on the IEEE 39-bus test system using PowerWorld simulator.}

\end{abstract}

\section{Introduction}
\label{Introduction}
\mrn{In recent years, the use of communication and control devices in electric power grids has significantly improved power systems' operation, enhancing both reliability, controllability, and efficiency. However, these systems have also become more vulnerable due to the presence of cyber elements with a high potential for anomaly intrusion. False Data Injection (FDI) attacks have been recognized as an important category of cyberattack that threatens the secure operation of power systems. An FDI attack is a type of cyberattack designed to deliver false information, intending to mislead or disrupt system operations~\cite{tran2021design}. In a general FDI attack scenario, an attacker injects an arbitrary amount of error into the system's state estimator by hacking sensors and measurement units~\cite{iranpour2024ac, iranpour2024designing}. If this attack vector fulfill certain conditions, the FDI will pass commonly used residue-based bad data detectors while evades system estimator and misleading it~\cite{liu2011false,boyaci2021joint,boyaci2022generating, boyaci2022infinite, boyaci2022spatio}.}

\mrn{Several studies over the past decade have reviewed FDI attacks from various perspectives. In addition to construction methods, goals, and consequences of FDIAs from the attackers' perspective, as well as detection and defense strategies reviewed in~\cite{zhang2019false}, FDI attacks can be analyzed based on their impacts on power system operations~\cite{chung2018local}. Depending on various factors in designing an FDI attack, such as locations, numbers, sequence of injecting false data, or values of the compromised devices, the impact of FDI attacks may range from vandalizing critical infrastructure to gaining specific economic advantages or benefits~\cite{tran2021design}. Generally, after conducting an FDI attack, the state estimator generates a wrong set of inputs for various control algorithms, such as optimal power flows, contingency analysis, economic dispatch, stabilizer units, or other controllers in the power system. Subsequently, these incorrectly generated control actions might initiate the system destabilization process, possibly leading to a collapse~\cite{xu2020review}. }

\mrn{Power flow algorithms are essential for the reliable operation of power systems~\cite{narimani2023tightening,narimani2020tightening, narimani2020strengthening,narimani2018comparison, narimani2018empirical,narimani2018improving}. Attackers can mislead operators by adjusting some power flow measurements, making them appear as real measurements. Based on this incorrect analysis, the controlling algorithm may trigger the relay protection system to open certain breakers, isolating branches and interrupting the normal operation of the system, potentially leading to a power outage~\cite{chung2018local}. Numerous studies have investigated conducting FDI attacks on power systems. For example, in~\cite{chung2018local}, attackers designed an attack strategy based on AC power flow equations by selecting the target power transmission line. In~\cite{soltan2019line} and~\cite{davis2012power}, it has been shown that FDI attacks can disconnect multiple power transmission lines within the attack zone or cause significant changes in the power flow of the grid. Suspicious changes in power flow variations in power transmission lines due to FDI attacks, based on the Power Transfer Distribution Factor (PTDF), have been analyzed in~\cite{li2019enhancing}. The impacts of FDI attacks on the security, operation, and economics of power systems have been analyzed using the optimal power flow in~\cite{rahman2014formal}. It has been shown that if metering devices contain massively corrupted measurements, existing FDI attacks can be detected in the system~\cite{anwar2017modeling}, but these studies have not shown how FDI attacks can cause large variations in transmission line power flow. In all of these studies, the impacts of an FDI attack on transmission line power flow have not been empirically analyzed. A mathematical method based on matrix analysis and graph theory is presented in~\cite{mohammadi2022impact} to quantify how attacking a predetermined set of power transmission lines with a DC-based FDI attack can affect power flow changes.}

\mrn{Although~\cite{mohammadi2022impact} presented a mathematical analysis for investigating the impacts of an FDI attack on the transmission lines' power flow, their approach is based on the DC power flow, making their FDI attack easily detectable by BDD algorithms. To address this, we proposed an approach for designing an AC-based FDI attack in the power system. The impact of various FDI attacks differs on the power system; we recently designed an optimal AC FDI attack scheme to design FDI attacks with minimum changes in state variables~\cite{iranpour2024ac}. Designing an optimal FDI attack results in a lower residual estimation vector compared to the residuals of estimation with arbitrary attack vectors. Thus, the optimal FDI attacks have a higher chance of dodging the BDD algorithm. Although slight changes in state variables lead to lower residuals, the impact of the designed FDI attack on the power system is lesser. In this connection, the attacker must trade-off between having residuals and FDI attack impact. This paper analyzes the impacts of different AC FDI attacks on the transmission lines' power flow. We used PowerWorld software to demonstrate the impact of FDI attacks on the power flows in power transmission lines.}


\mrn{The rest of this paper is organized as follows: Section~\ref{sec:design_AC_FDI_ATTACK} briefly describes the approach for designing AC FDI attacks and how the designed attacks impact power systems.
Section~\ref{sec:ANALYSIS}, presents the simulation results of the IEEE 39-bus test system. Section~\ref{sec:conclusion} concludes the paper.}

\section{AC FDI attack Impact on Power Systems}
\label{sec:design_AC_FDI_ATTACK}
\mrn{Analyzing the impacts of AC False Data Injection (FDI) attacks on power systems involves a systematic process to understand how such attacks affect system operation. The first step in this process is to understand the mechanism of designing an AC FDI attack, identifying the components and corresponding variables that need to be manipulated for a successful attack. In the second step, based on the injected attack vector to the attacked measurements, the values of these measurements and the corresponding estimation residuals of these variables (after estimating based on the manipulated measurements) will change. These attack vectors and estimation residual values can vary based on the magnitude of the attack vector's components, i.e., how severe the FDI attacks are. If these injected FDI attacks are optimally designed, the corresponding residuals will be less than the residuals for random FDI attacks. In this context, designing an optimal AC FDI attack vector is introduced in the following section. Finally, the proposed attack vector design will be implemented in the PowerWorld simulator to analyze the impacts of FDI attacks on the power flows of transmission lines as an index that can reflect the impacts of an attack on power systems. }

\subsection{Designing AC FDI Attack}

\mrn{The Remote Terminal Units (RTUs) and Phasor Measurement Units (PMUs) collect various measurements, including power flows, power injections, voltage magnitudes, voltage angles, and current magnitudes~\cite{wood2013power}. These measurements in the power systems are transmitted to the control center for state estimation. State Estimation (SE) is an important process in control centers to infer the status of various power system state variables (including voltage magnitudes and angles at buses) to deduce the status of the power system in different operation steps.
This metering process is prone to errors, which could affect the SE operation. To detect these errors, the difference between the observed measurement ($z$) and the estimated measurements ($h(x)$) is considered as an index, referred to as the residuals of estimations. If these residuals are greater than a predefined threshold ($r=|z - h(x)|>\tau$), it indicates the presence of bad data in the system. However, attackers can inject false data, as presented in equation~\eqref{eq:AC_residual}, capable of bypassing the BDD algorithms.}

\vspace{-.4cm}
{\begin{align}
\label{eq:AC_residual}
&\nonumber r_a = z_a-h(x_a) \\
&\nonumber ~~~= z_a- h(x_a)+h(x)-h(x)\\
&~~~\nonumber= z+a-h(x_a)+h(x)-h(x)\\ &~~~= r+a- h(x_a)+h(x).
\end{align}}

\mrn{In this equation, $r_a$, $r$, $x_a$ and $x$ represent the residuals for under-attack and normal measurements, and state vectors for under-attack and normal measurement vectors, respectively. Additionally, $h$ denotes the set of nonlinear power flow equations that relate measurements and states. In this scenario, if vector $a$ is defined as in equation~\eqref{eq:AC_undetectable}, this crafted false data injection vector can disrupt the data distribution while evading the BDD mechanism.}

\vspace{-.4cm}
{\begin{align}
\label{eq:AC_undetectable}
    a=h(x_a)-h(x).
\end{align}}
\vspace{-.4cm}

\mrn{It is notable that in equation~\eqref{eq:AC_undetectable}, $x_a$ is equal to $x+c$, where vector $c$ is an arbitrary $n\times 1$ non-zero vector and $n$ is the number of measurements that needs to be changed to conduct the attack.}
\mrn{Defining the nonlinear function $h(x)$ and designing the attack vector  $a$ include some assumptions and steps which are listed below:}
\vspace{-.1cm}
\mrn{
\begin{itemize}
\item In designing a successful AC FDI attack, the main point is rationalizing any changes in power flow equations that stem from injecting an attack vector.
\item Attacked and intact zones should be determined, and any changes in the attacked area and between the attacked area and the intact area should be rationalized. 
\item The attack area includes the region directly targeted by the attack, consisting of a set of buses, an enclosed zone where all changes resulting from the attack occur internally. Meanwhile, the intact area remains unaffected by the FDI attack.
\item In this scenario, when executing an FDI attack on a targeted region, it is crucial to maintain the total power transfers between affected and unaffected zones equivalent before and after the attack scenarios. To meet this requirement, the attack area should be surrounded by buses with nonzero power injection. This accounts for every power change exclusively within the attack zone by crafting injection measurements on these buses~\cite{hug2012vulnerability}. Additionally, to ensure that all changes remain within the attack area during the attack, the state variables at the boundary must remain unchanged. Therefore, after designating specific buses or groups of buses as primary focal points to determine the attack zone, this focal bus will be expanded by including zero injection buses and considering neighboring buses with non-zero injections as the boundary of the attack zone.
\item Within the attacked zone, the algebraic sum of generated and consumed power must remain unchanged. Therefore, in zero-injection buses within the attack zone, the algebraic sum of active/reactive power flows after the attack must be zero (complex power flows must be zero). In non-zero injection buses within the attack zone, the power injection after the attack equals the primary injection power plus the sum of all changes in power flows of lines connected to this non-zero injection bus that exists within the attack zone. This assumption has been presented in equations~\ref{eq:p_non_zero_injection} and \ref{eq:q_non_zero_injection}.
\end{itemize}}

\vspace{-.4cm}
\begin{subequations}
\begin{small}
\begin{align}
\label{eq:p_non_zero_injection}
&\tilde{P}_{m}=P_{m}+\sum_{(m,l)\in \mathcal{L}_A}(\tilde{P}_{m,l}-P_{m,l}),\\
\label{eq:q_non_zero_injection}
&\tilde{Q}_{m}=Q_{m}+\sum_{(m.l)\in \mathcal{L_A}}(\tilde{Q}_{m,l}-Q_{m,l}).
\end{align}
\end{small}
\end{subequations}

\mrn{In these Equations, $\tilde{P_{m}}$, $\tilde{Q_{m}}$, $P_{m}$, and $Q_{m}$ represent the active and reactive power injections at bus $m$ after and before the attack, respectively. Also, $\tilde{P}{m,l}$ and $\tilde{Q}{m,l}\in \mathcal{L}_A$ are the active and reactive powers of the lines within the attack zone that are connected to bus $m$.  $\mathcal{L}_A$ is the set of lines that are inside the attack zone. }

\subsection{Comparing the Impact of Various FDI attacks}

\mrn{Various factors can be considered to analyze the impacts of different FDI attacks on power systems. We have focused on two factors: the attack zones and manipulated measurements within those zones. Our goal is to examine how these factors affect the state estimation modules when false data is injected. Specifically, we assess the impact of FDI attacks on power flow through transmission lines by monitoring changes in line capacity after injecting FDI attack vectors into corresponding measurements within the attack zone. Our aim is to demonstrate the effects on transmission line power flow after manipulating state estimation with false data and subsequently adjusting power flows at the control center based on this information. The residuals of the estimated variables are considered as the second factor for analyzing the impact of FDI attacks. The magnitudes of the components of the attack vector are a crucial factor that can affect the residuals of the estimations after the attack~\cite{iranpour2024ac}. As described in~\cite{iranpour2024ac}, designing an optimal attack vector results in residuals of corresponding variables in the state estimation module being less than those from an arbitrary AC FDI attack. Further details on optimal design for AC FDI attacks can be found in~\cite{iranpour2024ac}. }

\subsection{Implementing the FDI attack in the PowerWorld Simulator}
\mrn{This paper utilizes the PowerWorld simulator to assess the impacts of optimal and arbitrary AC FDI attacks on transmission line power flow. To do so, we inject the designed attack vector into the state variables in the PowerWorld simulator. A critical aspect of successfully implementing this vector involves fixing the voltages of the boundary buses in the attack zones within PowerWorld. Subsequently, the voltage (magnitude and angle) attack vector is injected into the buses within the attack zone. 
Note that to fix the voltage magnitude at the boundary buses in the PowerWorld Simulator, we utilized voltage regulator devices. These devices help maintain the desired voltage levels at the boundary buses, ensuring that the power flow into and out of the attack zone remains unchanged. This step is essential for accurately simulating and analyzing the impacts of AC FDI attacks on the power system. The analysis of transmission line power flow for optimal and arbitrary AC FDI attacks is based on the occupied capacity of the lines in PowerWorld. Further details are presented in Section~\ref{sec:ANALYSIS}.  }

\section{Simulation Results and Discussion}
\label{sec:ANALYSIS}
\mrn{In this section, we demonstrate the impacts of FDI attacks on transmission line power flows by implementing optimal and arbitrary AC FDI attacks on the IEEE 39-bus test case from the PGLib-OPF v18.08 benchmark library~\cite{pglib}. Initially, we design optimal and arbitrary AC FDI attack vectors based on the approach described in Section \ref{sec:design_AC_FDI_ATTACK}. Subsequently, we implement these designed attack vectors in the PowerWorld simulator to assess the impacts of AC FDI attacks. The steps required for implementing FDI attacks and their corresponding impacts on power systems are discussed in the following.}


\mrn{The first step in conducting an AC FDI attack in the proposed approach is identifying the attack zone. To this end, we considered buses $3,18,17,16,15,24,21,27,26,28,29,25$ as the attack zone, with buses $3, 16, 15, 21, 24, 25$, and $29$ as the boundary buses. Five buses, including $17$, $18$, $26$, $27$, and $28$, are within the attack zone and thus have variable voltage values. These values, required for conducting an optimal AC FDI attack, are calculated based on an optimization problem with constraints demonstrated in Section \ref{sec:design_AC_FDI_ATTACK}. More details about designing optimal AC FDI attacks in power systems can be found in our recent work~\cite{iranpour2024ac}.
We design the arbitrary AC FDI attack by fixing the objective function in the optimization problem for designing the AC FDI attack. By doing this, all constraints for designing the AC FDI attack are respected. After designing the attacks, we can calculate the attack vector $a$ for voltage, power flow, and power injection measurements in the attack zone.}

\mrn{An additional constraint for overloading the line between buses $26$ and $27$ is added to the AC FDI attack design problem, which enforces the power flow in the line to be greater than or equal to $1.3$ times its power flow before conducting the attack. The value of the active power flow of this line is equal to $2.573 (p.u.)$, $3.3466 (p.u.)$, and $11.4067 (p.u.)$ for the before attack, optimal attack, and arbitrary attack situations, respectively.}

\begin{figure}
    \centering
\includegraphics[scale=0.454]{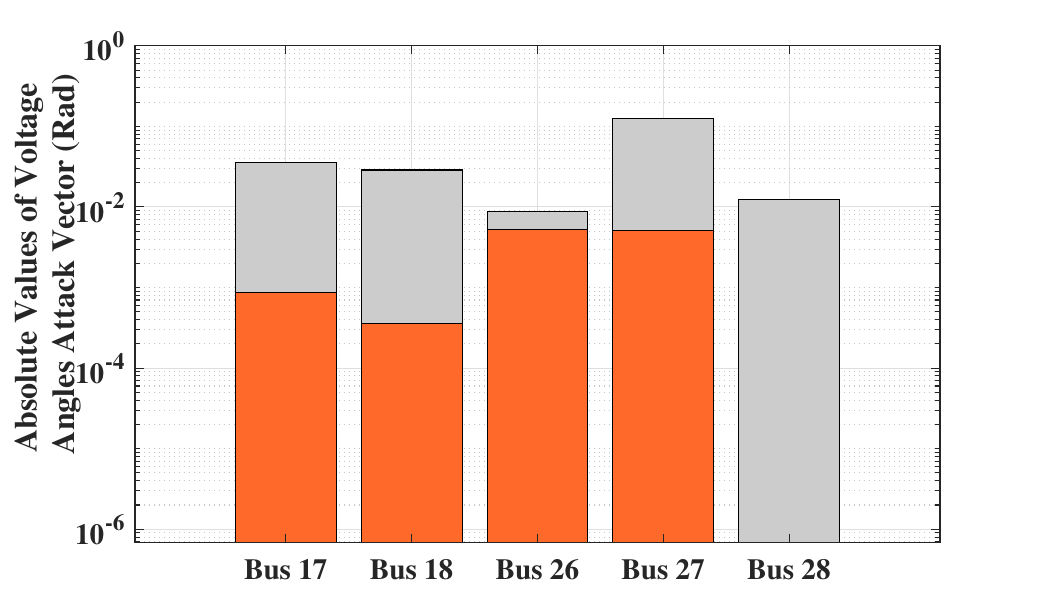}
\vspace{-.4cm}
	\caption{\mrn{Attack Vector for voltage angle. The red and gray bars indicate the voltage angles that need to be added to the measurement for optimal and
arbitrary AC FDI attacks, respectively}}
	\label{fig:voltage_ang39}
\end{figure}
 
\begin{figure}
    \centering
\includegraphics[scale=0.42]{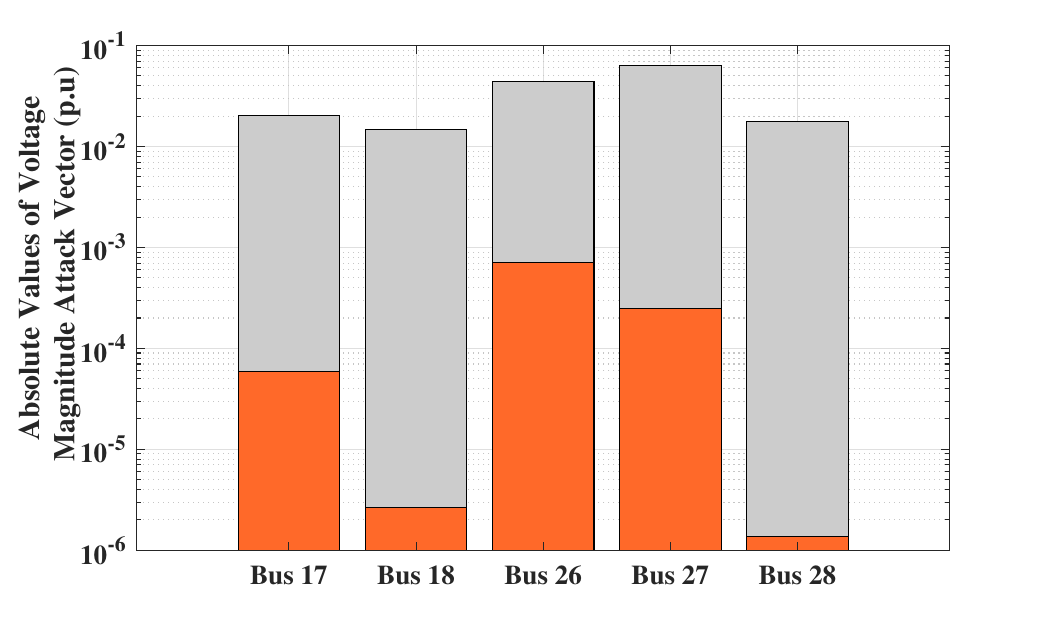}
\vspace{-.4cm}
	\caption{Attack Vector for voltage magnitude. The red and gray bars indicate the voltage magnitudes that need to be added to the measurement for optimal
and arbitrary FDI attacks, respectively}
	\label{fig:voltage_mag39}
\end{figure}

\mrn{The corresponding attack vectors for voltage magnitude and angle for the optimal and arbitrary AC FDI attacks are shown in Figs.~\ref{fig:voltage_ang39} and~\ref{fig:voltage_mag39}, respectively. From these figures, it is clear that the magnitude of the attack vector for both voltage magnitude and angle in the optimal AC FDI attack is smaller than their corresponding values for the arbitrary AC FDI attack. The objective function for designing the AC FDI attack in our recent work in \cite{iranpour2024ac} is to minimize the difference between state variables before and after conducting the AC FDI attack. This, causes the attack vector to be as small as possible while satisfying the constraints, including the overloading constraint that enforces the power flow in the lines to be greater than a predefined value.}

\begin{table*}
\caption{Voltage magnitude and angle of busses in the attack zone for before attack, optimal attack and arbitrary attack situations}
\label{Table:Voltages}
\begin{tabular}{|c|c|c|c|c|c|c|}
\hline 
Bus Number & \multicolumn{2}{c|}{Before attack} & \multicolumn{2}{c|}{Optimal attack} & \multicolumn{2}{c|}{Arbitrary attack}\tabularnewline
\hline 
\hline 
 & Voltage Magnitude & Voltage Angle & Voltage Magnitude & Voltage Angle & Voltage Magnitude & Voltage Angle\tabularnewline
\hline 
{\cellcolor{shadecolor}} 3 &{\cellcolor{shadecolor}} 1.0307 &{\cellcolor{shadecolor}} -12.2763 &{\cellcolor{shadecolor}} 1.0307 &{\cellcolor{shadecolor}} -12.2764 &{\cellcolor{shadecolor}} 1.0307 &{\cellcolor{shadecolor}} -12.2764\tabularnewline
\hline 
{\cellcolor{shadecolor}} 15 &{\cellcolor{shadecolor}} 1.0161 &{\cellcolor{shadecolor}} -11.3453 &{\cellcolor{shadecolor}} 1.0161 &{\cellcolor{shadecolor}} -11.3454 &{\cellcolor{shadecolor}} 1.0162 &{\cellcolor{shadecolor}} -11.3454\tabularnewline
\hline 
{\cellcolor{shadecolor}} 16 &{\cellcolor{shadecolor}} 1.0325 &{\cellcolor{shadecolor}} -10.0333 &{\cellcolor{shadecolor}} 1.0325 &{\cellcolor{shadecolor}} -10.0333 &{\cellcolor{shadecolor}} 1.0325 &{\cellcolor{shadecolor}} -10.0333\tabularnewline
\hline 
17 & 1.0342 & -11.1164 & 1.0342 & -11.1659 & 1.014 & -13.1662\tabularnewline
\hline 
18 & 1.0315 & -11.9861 & 1.0316 & -11.9655 & 1.0167 & -13.634\tabularnewline
\hline 
{\cellcolor{shadecolor}} 21 &{\cellcolor{shadecolor}} 1.0323 &{\cellcolor{shadecolor}} -7.6287 &{\cellcolor{shadecolor}} 1.0323 &{\cellcolor{shadecolor}} -7.6287 &{\cellcolor{shadecolor}} 1.0323 &{\cellcolor{shadecolor}} -7.6287\tabularnewline
\hline 
{\cellcolor{shadecolor}} 24 &{\cellcolor{shadecolor}} 1.0380 &{\cellcolor{shadecolor}} -9.9137 &{\cellcolor{shadecolor}} 1.0380 &{\cellcolor{shadecolor}} -9.9138 &{\cellcolor{shadecolor}} 1.0380 &{\cellcolor{shadecolor}} -9.9138\tabularnewline
\hline 
{\cellcolor{shadecolor}} 25 &{\cellcolor{shadecolor}} 1.0576 &{\cellcolor{shadecolor}} -8.3692 &{\cellcolor{shadecolor}} 1.0576 &{\cellcolor{shadecolor}} -8.3692 &{\cellcolor{shadecolor}} 1.0577 & {\cellcolor{shadecolor}} -8.3692\tabularnewline
\hline 
26 & 1.0525 & -9.4387 & 1.0533 & -9.137 & 1.0081 & -8.9311\tabularnewline
\hline 
27 & 1.0383 & -11.3621 & 1.0381 & -11.6541 & 0.9749 & -18.5778\tabularnewline
\hline 
28 & 1.0503 & -5.9283 & 1.0504 & -5.9284 & 1.0328 & -5.2174\tabularnewline
\hline 
{\cellcolor{shadecolor}} 29 &{\cellcolor{shadecolor}} 1.0501 &{\cellcolor{shadecolor}} -3.1698 &{\cellcolor{shadecolor}} 1.0501 &{\cellcolor{shadecolor}} -3.1699 &{\cellcolor{shadecolor}} 1.0501 &{\cellcolor{shadecolor}} -3.1699\tabularnewline
\hline 
\end{tabular}
\end{table*}


\mrn{The voltage magnitude and angle for buses in the attack zones, including those within the attack zone and boundary buses, are tabulated in Table~\ref{Table:Voltages}. From this table, it is clear that the voltage magnitude and angle for the boundary buses remain unchanged after the attack. The corresponding rows for the boundary buses are highlighted in this table, demonstrating that their values are the same before and after conducting the AC FDI attack. This is crucial for conducting successful AC FDI attacks, as the net power flow into the attack zone must not change. Similarly, the power injections for the non-zero injection buses within the attack zone are tabulated in Table~\ref{Table:Power_Injection}. The active and reactive power flows in the lines are shown in Table~\ref{Table:Power_flow} for selected branches. The power flow of the branches connected to the boundary buses that are outside the attack zone should not change after conducting the AC FDI attack to ensure the net zero power injection into the attack zone. These branches are highlighted in Table~\ref{Table:Power_flow}. It is evident that the proposed algorithm successfully handles this condition and does not change the flow of these lines during the process of designing the AC FDI attack.}

\begin{figure}
 \centering
    \captionsetup{justification=centering}
\begin{tikzpicture}
\node at (0,0) {\includegraphics[scale=0.8,trim= 8cm 5.6cm 3.3cm 5.4cm  2.3cm,clip]{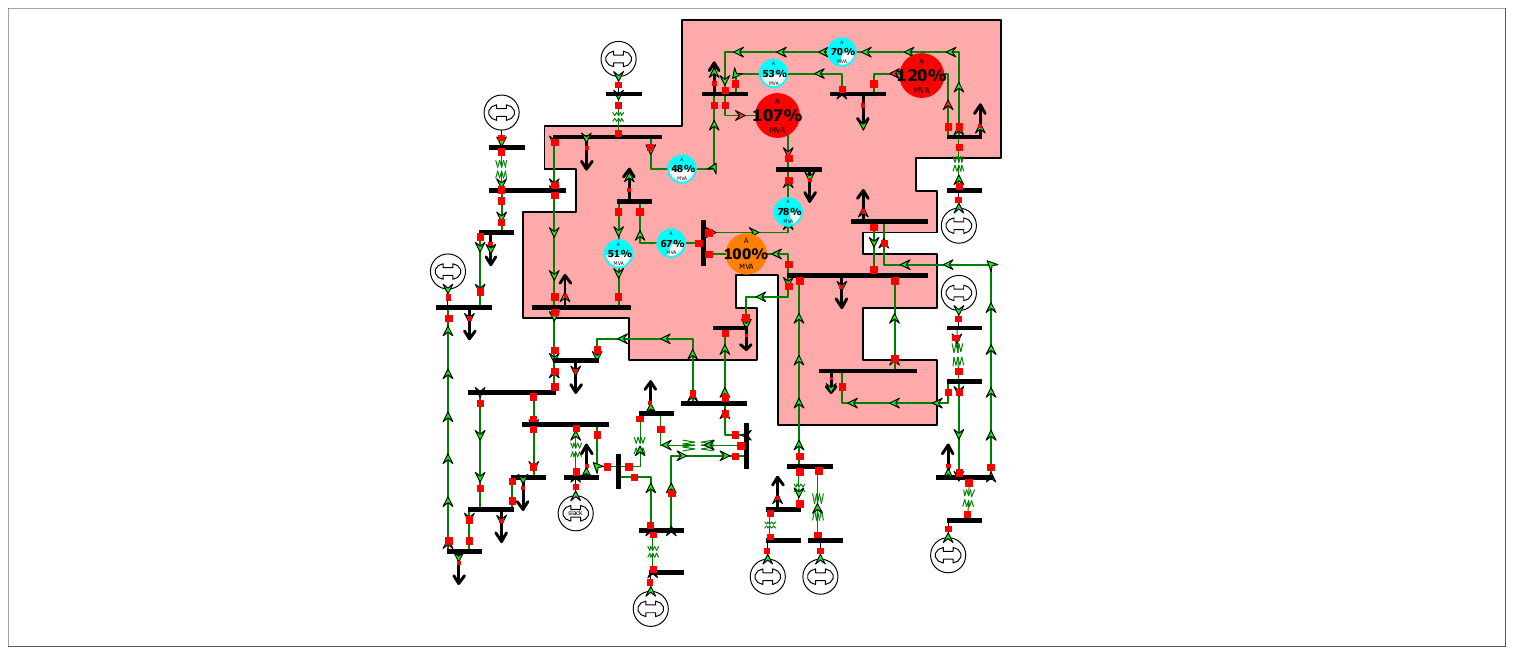}};

\node [anchor=west, font=\scriptsize] at (-1.4,-1.7) {$B_{19}$};
\node [anchor=west, font=\scriptsize] at (-0.8,-3) {$B_{33}$};
\node [anchor=west, font=\scriptsize] at (-4.6,3.3) {$B_{37}$};
\node [anchor=west, font=\scriptsize] at (-6.2,2.3) {$B_{30}$};
\node [anchor=west, font=\scriptsize] at (-6.1,1.8) {$B_{2}$};
\node [anchor=west, font=\scriptsize] at (-6.2,1.3) {$B_{1}$};
\node [anchor=west, font=\scriptsize] at (-6.9,0.3) {$B_{39}$};
\node [anchor=west, font=\scriptsize] at (-4.4,0.4) {$B_{3}$};
\node [anchor=west, font=\scriptsize] at (-4.4,1.6) {$B_{18}$};
\node [anchor=west, font=\scriptsize] at (-3.1,1.6) {$B_{17}$};
\node [anchor=west, font=\scriptsize] at (-4.2,2.3) {$B_{25}$};
\node [anchor=west, font=\scriptsize] at (-3.2,3.3) {$B_{26}$};
\node [anchor=west, font=\scriptsize] at (-1.2,2.9) {$B_{28}$};
\node [anchor=west, font=\scriptsize] at (0,2.5) {$B_{29}$};
\node [anchor=west, font=\scriptsize] at (1,1.9) {$B_{38}$};
\node [anchor=west, font=\scriptsize] at (-1.5,2.3) {$B_{27}$};
\node [anchor=west, font=\scriptsize] at (-0.4,1.6) {$B_{24}$};
\node [anchor=west, font=\scriptsize] at (-1.2,0.85) {$B_{16}$};
\node [anchor=west, font=\scriptsize] at (-1.2,-0.4) {$B_{21}$};
\node [anchor=west, font=\scriptsize] at (-2.8,0.2) {$B_{15}$};
\node [anchor=west, font=\scriptsize] at (-2.9,-0.8) {$B_{14}$};
\node [anchor=west, font=\scriptsize] at (-4.6,-0.3) {$B_{4}$};
\node [anchor=west, font=\scriptsize] at (-5.8,-0.8) {$B_{5}$};
\node [anchor=west, font=\scriptsize] at (-4.8,-1.2) {$B_{6}$};
\node [anchor=west, font=\scriptsize] at (-4,-1) {$B_{12}$};
\node [anchor=west, font=\scriptsize] at (-2.1,-1.5) {$B_{13}$};
\node [anchor=west, font=\scriptsize] at (-2.4,-3) {$B_{34}$};
\node [anchor=west, font=\scriptsize] at (-2.4,-2.5) {$B_{20}$};
\node [anchor=west, font=\scriptsize] at (-3.4,-3) {$B_{10}$};
\node [anchor=west, font=\scriptsize] at (-3.2,-3.6) {$B_{32}$};
\node [anchor=west, font=\scriptsize] at (-6.5,-3.3) {$B_{9}$};
\node [anchor=west, font=\scriptsize] at (-5.8,-2.3) {$B_{8}$};
\node [anchor=west, font=\scriptsize] at (-5.2,-2.3) {$B_{7}$};
\node [anchor=west, font=\scriptsize] at (-4.9,-1.8) {$B_{31}$};
\node [anchor=west, font=\scriptsize] at (-4.2,-2.3) {$B_{11}$};
\node [anchor=west, font=\scriptsize] at (0,0) {$B_{35}$};
\node [anchor=west, font=\scriptsize] at (1,-0.8) {$B_{22}$};
\node [anchor=west, font=\scriptsize] at (-0.2,-2) {$B_{23}$};
\node [anchor=west, font=\scriptsize] at (,-2.8) {$B_{36}$};
\end{tikzpicture}%
\vspace{-5pt}
 \caption{Percentage of the Occupied
Capacity of the Lines based on the Power Flows of Lines in arbitrary attack scenario}
\label{fig:powerworld_nonoptimum}
\end{figure}
\begin{figure}
    \centering
\includegraphics[scale=0.454]{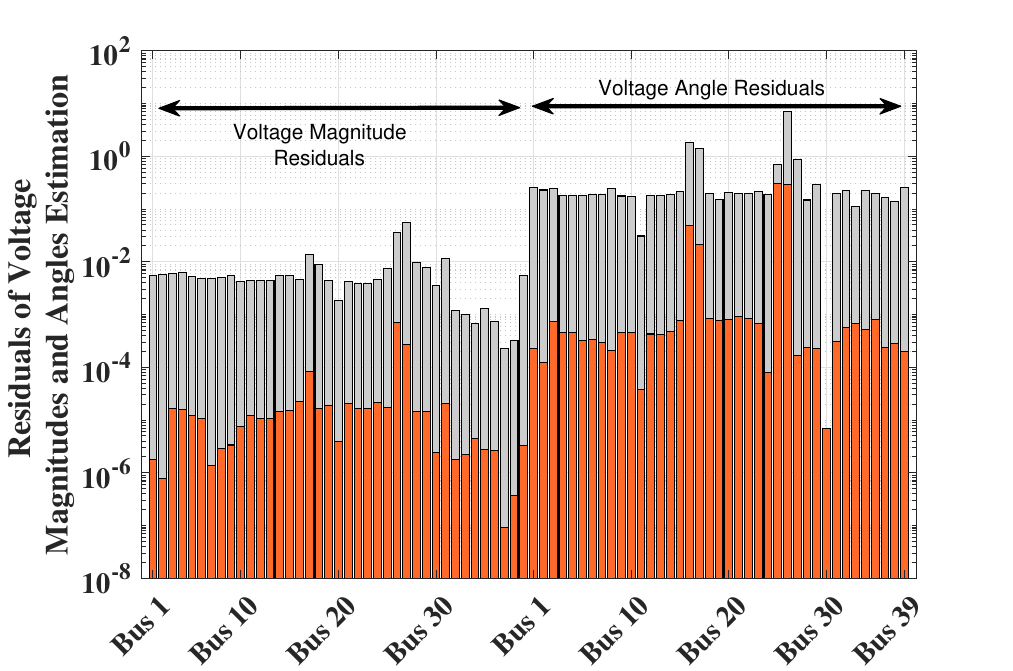}
	\caption{Residuals of AC state estimator for voltage magnitudes and angles. The red
and gray bars indicate the residuals of the AC state estimator for voltage
magnitudes under optimal and arbitrary FDI attacks, respectively.}
	\label{fig:residual_ang39}
\end{figure}
\begin{figure}
 \centering
    \captionsetup{justification=centering}
\begin{tikzpicture}
\node at (0,0) {\includegraphics[scale=0.8,trim= 8cm 5.4cm 3cm 5.5cm  2.3cm,clip]{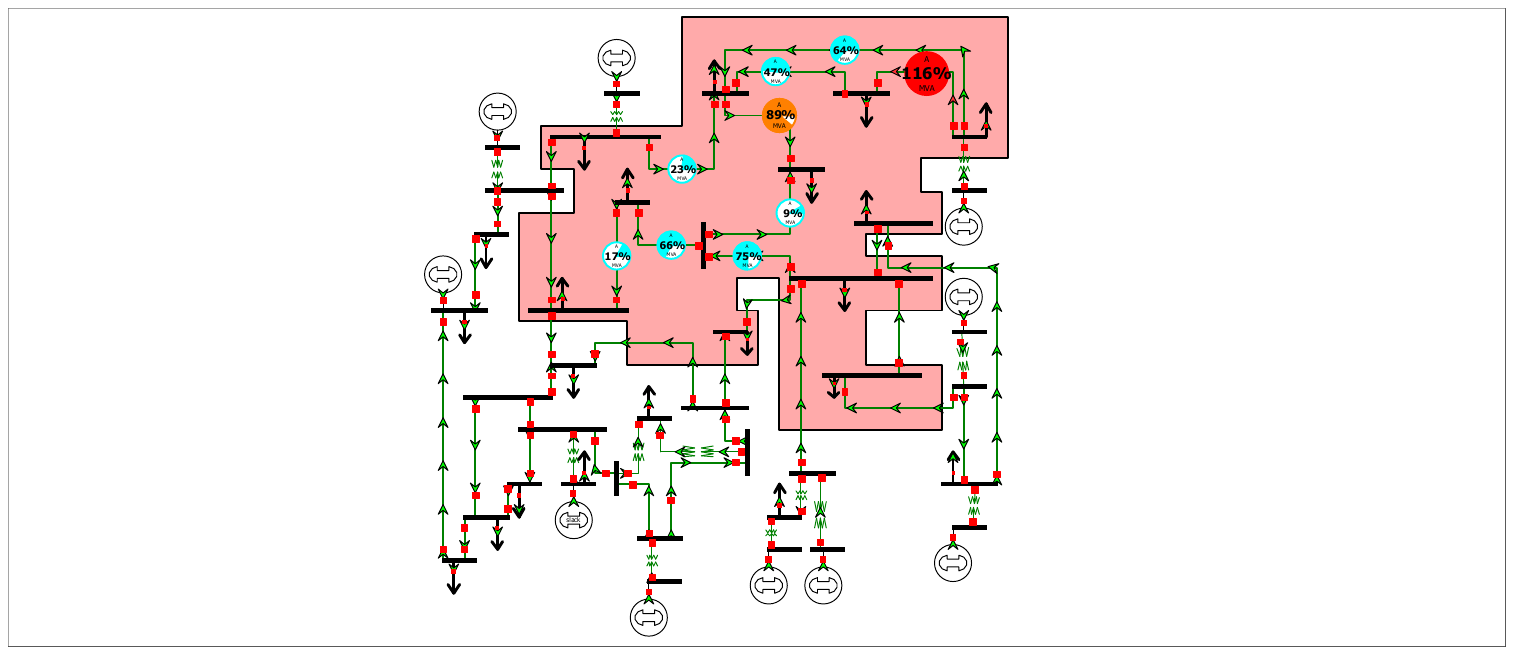}};

\node [anchor=west, font=\scriptsize] at (-1.4,-1.8) {$B_{19}$};
\node [anchor=west, font=\scriptsize] at (-0.8,-3) {$B_{33}$};
\node [anchor=west, font=\scriptsize] at (-4.7,3.3) {$B_{37}$};
\node [anchor=west, font=\scriptsize] at (-6.2,2.3) {$B_{30}$};
\node [anchor=west, font=\scriptsize] at (-6.1,1.8) {$B_{2}$};
\node [anchor=west, font=\scriptsize] at (-6.4,1.3) {$B_{1}$};
\node [anchor=west, font=\scriptsize] at (-7,0.3) {$B_{39}$};
\node [anchor=west, font=\scriptsize] at (-4.5,0.4) {$B_{3}$};
\node [anchor=west, font=\scriptsize] at (-4.5,1.6) {$B_{18}$};
\node [anchor=west, font=\scriptsize] at (-3.1,1.6) {$B_{17}$};
\node [anchor=west, font=\scriptsize] at (-4.2,2.4) {$B_{25}$};
\node [anchor=west, font=\scriptsize] at (-3.35,3.3) {$B_{26}$};
\node [anchor=west, font=\scriptsize] at (-1.2,2.9) {$B_{28}$};
\node [anchor=west, font=\scriptsize] at (0,2.5) {$B_{29}$};
\node [anchor=west, font=\scriptsize] at (1,1.9) {$B_{38}$};
\node [anchor=west, font=\scriptsize] at (-1.5,2.35) {$B_{27}$};
\node [anchor=west, font=\scriptsize] at (-0.5,1.7) {$B_{24}$};
\node [anchor=west, font=\scriptsize] at (-1.2,0.85) {$B_{16}$};
\node [anchor=west, font=\scriptsize] at (-1.2,-0.4) {$B_{21}$};
\node [anchor=west, font=\scriptsize] at (-2.8,0.2) {$B_{15}$};
\node [anchor=west, font=\scriptsize] at (-3,-0.78) {$B_{14}$};
\node [anchor=west, font=\scriptsize] at (-4.8,-0.3) {$B_{4}$};
\node [anchor=west, font=\scriptsize] at (-5.8,-0.7) {$B_{5}$};
\node [anchor=west, font=\scriptsize] at (-4.8,-1.2) {$B_{6}$};
\node [anchor=west, font=\scriptsize] at (-4.1,-1) {$B_{12}$};
\node [anchor=west, font=\scriptsize] at (-2.2,-1.6) {$B_{13}$};
\node [anchor=west, font=\scriptsize] at (-2.5,-3) {$B_{34}$};
\node [anchor=west, font=\scriptsize] at (-2.5,-2.5) {$B_{20}$};
\node [anchor=west, font=\scriptsize] at (-3.4,-3) {$B_{10}$};
\node [anchor=west, font=\scriptsize] at (-3.2,-3.6) {$B_{32}$};
\node [anchor=west, font=\scriptsize] at (-6.7,-3.3) {$B_{9}$};
\node [anchor=west, font=\scriptsize] at (-6,-2.75) {$B_{8}$};
\node [anchor=west, font=\scriptsize] at (-5.2,-2.3) {$B_{7}$};
\node [anchor=west, font=\scriptsize] at (-4.5,-2.3) {$B_{31}$};
\node [anchor=west, font=\scriptsize] at (-4.2,-1.7) {$B_{11}$};
\node [anchor=west, font=\scriptsize] at (0,0) {$B_{35}$};
\node [anchor=west, font=\scriptsize] at (0,-0.8) {$B_{22}$};
\node [anchor=west, font=\scriptsize] at (-0.2,-2) {$B_{23}$};
\node [anchor=west, font=\scriptsize] at (,-2.8) {$B_{36}$};

\end{tikzpicture}%
\vspace{-.4cm}
 \caption{Percentage of the Occupied
Capacity of the Lines based on the Power Flows of Lines in optimal attack scenario}
	\label{fig:powerworld_optimum}
\end{figure}

\mrn{Note that, in these tables, the magnitude of attack values for the arbitrary attacks is greater than the corresponding values in the optimal attack scenarios. The impacts of these differences are observable in the residuals of the corresponding estimated parameters. Residuals of voltage magnitude and angle estimations for optimal and arbitrary scenarios are shown in Fig.~\ref{fig:residual_ang39} for both optimal and arbitrary attacks. These residuals are calculated by first injecting the attack vector into the corresponding measurements and then using the MATGRID toolbox~\cite{cosovic2024matgrid} to perform AC state estimation. Note that the difference between the residuals of arbitrary and optimal attacks is noticeable; the residuals for the arbitrary attack are much larger than those for the optimal attack, making the arbitrary attack an easy target to be detected by the BDD algorithms. The impact of the arbitrary and optimal attacks on the power flows is shown in Figs.~\ref{fig:powerworld_nonoptimum} and \ref{fig:powerworld_optimum}, respectively. We want to emphasize that although the impact of the arbitrary attack on the power flow in the grid is larger than that of the optimal attack, their residuals are also greater than those for the optimal attack.}



\vspace{-.4cm}

\begin{table*}
\caption{Active and reactive power injection in no-zero injecting busses in the attack zone for before attack, optimal attack and arbitrary attack situations}
\label{Table:Power_Injection}
\begin{tabular}{|c|c|c|c|c|c|c|}
\hline 
\multirow{3}{*}{\thead{Bus Number}} & \multicolumn{2}{c|}{\thead{Before Attack}} & \multicolumn{2}{c|}{\thead{Optimal Attack}} & \multicolumn{2}{c|}{\thead{Arbitrary Attack}}\tabularnewline
\cline{2-7} \cline{3-7} \cline{4-7} \cline{5-7} \cline{6-7} \cline{7-7} 
 & \multirow{2}{*}{\thead{Active Power}} & \multirow{2}{*}{\thead{Reactive Power}} & \multirow{2}{*}{\thead{Active Power}} & \multirow{2}{*}{\thead{Reactive Power}} & \multirow{2}{*}{\thead{Active Power}} & \multirow{2}{*}{\thead{Reactive Power}}\tabularnewline
 & \thead{Injection (p.u)} &  \thead{Injection (p.u)}& \thead{Injection (p.u)} & \thead{Injection (p.u)} & \thead{Injection (p.u)} & \thead{Injection (p.u)}\tabularnewline
\hline 
3 & -3.22 & -0.024 & -3.2487 & -0.0217 & -0.8675 & 0.9514\tabularnewline
\hline 
16 & -3.29 & -0.323 & -3.1866 & -0.3222 & 1.0409 & 1.8335\tabularnewline
\hline 
18 & -1.58 & -0.3 & -1.3926 & -0.3101 & -2.9305 & -0.639\tabularnewline
\hline 
25 & -2.24 & -0.472 & -2.422 & -0.4801 & -2.4141 & 0.9969\tabularnewline
\hline 
26 & -1.39 & -0.17 & -0.2226 & -0.1528 & 7.7404 & -1.4558\tabularnewline
\hline 
27 & -2.81 & -0.755 & -3.8398 & -0.7094 & -16.7259 & -1.9921\tabularnewline
\hline 
28 & -2.06 & -0.276 & -2.1829 & -0.2878 & -1.185 & -0.9875\tabularnewline
\hline 
29 & -2.835 & -0.269 & -2.9275 & -0.2822 & -3.8308 & 1.7312\tabularnewline
\hline 
\end{tabular}

\end{table*}

\begin{table*}
\caption{active and reactive power flow of lines in the attack zone for before attack,optimal attack and arbitrary attack situations}
\label{Table:Power_flow}
\begin{tabular}{|c|c|c|c|c|c|c|c|}
\hline 
\multicolumn{2}{|c|}{\thead{Branches}} & \multicolumn{2}{c|}{\thead{Before Attack}} & \multicolumn{2}{c|}{\thead{Optimal Attack}} & \multicolumn{2}{c|}{\thead{Arbitrary Attack}}\tabularnewline
\hline 
\hline 
\multirow{2}{*}{\thead{From Bus}} & \multirow{2}{*}{\thead{To Bus}} & \multirow{2}{*}{\thead{Active Power}} & \multirow{2}{*}{\thead{Reactive Power}} & \multirow{2}{*}{\thead{Active Power}} & \multirow{2}{*}{\thead{Reactive Power}} & \multirow{2}{*}{\thead{Active Power}} & \multirow{2}{*}{\thead{Reactive Power}}\tabularnewline
& & \thead{Flow (p.u)} &  \thead{Flow (p.u)}& \thead{Flow (p.u)} & \thead{Flow (p.u)} & \thead{Flow (p.u)} & \thead{Flow (p.u)}\tabularnewline
\hline 
{\cellcolor{shadecolor}} 2 &{\cellcolor{shadecolor}} 3 &{\cellcolor{shadecolor}} 3.1991 &{\cellcolor{shadecolor}} 0.8859 &{\cellcolor{shadecolor}} 3.1991 &{\cellcolor{shadecolor}} 0.8859 &{\cellcolor{shadecolor}} 3.1991 &{\cellcolor{shadecolor}} 0.8859\tabularnewline
\hline 
{\cellcolor{shadecolor}} 2 &{\cellcolor{shadecolor}} 25 &{\cellcolor{shadecolor}} -2.4459 &{\cellcolor{shadecolor}} 0.8297 &{\cellcolor{shadecolor}} -2.4459 &{\cellcolor{shadecolor}} 0.8297 &{\cellcolor{shadecolor}} -2.4459 &{\cellcolor{shadecolor}} 0.8297\tabularnewline
\hline 
{\cellcolor{shadecolor}} 3 & {\cellcolor{shadecolor}} 4 & {\cellcolor{shadecolor}}0.3734 & {\cellcolor{shadecolor}}1.1306 & {\cellcolor{shadecolor}}0.3734 & {\cellcolor{shadecolor}}1.1306 & {\cellcolor{shadecolor}}0.3734 & {\cellcolor{shadecolor}}1.1306\tabularnewline
\hline 
3 & 18 & -0.4076 & -0.1459 & -0.4363 & -0.1435 & 1.9449 & 0.8295\tabularnewline
\hline 
{\cellcolor{shadecolor}} 14 &{\cellcolor{shadecolor}} 15 &{\cellcolor{shadecolor}} 0.5031 &{\cellcolor{shadecolor}} -0.4068 &{\cellcolor{shadecolor}} 0.5031 &{\cellcolor{shadecolor}} -0.4068 &{\cellcolor{shadecolor}} 0.5031 &{\cellcolor{shadecolor}} -0.4068\tabularnewline
\hline 
{\cellcolor{shadecolor}} 15 &{\cellcolor{shadecolor}} 16 &{\cellcolor{shadecolor}} -2.6974 &{\cellcolor{shadecolor}} -1.5666 &{\cellcolor{shadecolor}} -2.6974 &{\cellcolor{shadecolor}} -1.5666 &{\cellcolor{shadecolor}} -2.6974 &{\cellcolor{shadecolor}} -1.5666\tabularnewline
\hline 
16 & 17 & 2.2402 & -0.4254 & 2.3436 & -0.4246 & 6.5711 & 1.7311\tabularnewline
\hline 
{\cellcolor{shadecolor}} 16 &{\cellcolor{shadecolor}} 19 &{\cellcolor{shadecolor}} -4.513 &{\cellcolor{shadecolor}} -0.542 &{\cellcolor{shadecolor}} -4.513 &{\cellcolor{shadecolor}} -0.542 &{\cellcolor{shadecolor}} -4.513 &{\cellcolor{shadecolor}} -0.542\tabularnewline
\hline 
{\cellcolor{shadecolor}} 16 &{\cellcolor{shadecolor}} 21 &{\cellcolor{shadecolor}} -3.296 &{\cellcolor{shadecolor}} 0.1444 &{\cellcolor{shadecolor}} -3.296 &{\cellcolor{shadecolor}} 0.1444 &{\cellcolor{shadecolor}} -3.296 & {\cellcolor{shadecolor}} 0.1444\tabularnewline
\hline 
{\cellcolor{shadecolor}} 16 &{\cellcolor{shadecolor}} 24 &{\cellcolor{shadecolor}} -0.4268 &{\cellcolor{shadecolor}} -0.9733 &{\cellcolor{shadecolor}} -0.4268 &{\cellcolor{shadecolor}} -0.9733 &{\cellcolor{shadecolor}} -0.4268 &{\cellcolor{shadecolor}} -0.9733\tabularnewline
\hline 
17 & 18 & 1.9904 & 0.1105 & 1.8313 & 0.1139 & 0.9912 & -0.4829\tabularnewline
\hline 
17 & 27 & 0.2464 & -0.4356 & 0.5086 & -0.4421 & 5.5494 & 1.9669\tabularnewline
\hline 
{\cellcolor{shadecolor}} 21 &{\cellcolor{shadecolor}} 22 &{\cellcolor{shadecolor}} -6.0442 &{\cellcolor{shadecolor}} -0.8726 &{\cellcolor{shadecolor}} -6.0442 &{\cellcolor{shadecolor}} -0.8726 &{\cellcolor{shadecolor}} -6.0442 &{\cellcolor{shadecolor}} -0.8726\tabularnewline
\hline 
{\cellcolor{shadecolor}} 23 &{\cellcolor{shadecolor}} 24 &{\cellcolor{shadecolor}} 3.5384 &{\cellcolor{shadecolor}} -0.005 &{\cellcolor{shadecolor}} 3.5384 &{\cellcolor{shadecolor}} -0.005 &{\cellcolor{shadecolor}} 3.5384 &{\cellcolor{shadecolor}} -0.005\tabularnewline
\hline 
25 & 26 & 0.6541 & -0.1881 & 0.4722 & -0.1962 & 0.48 & 1.2808\tabularnewline
\hline 
{\cellcolor{shadecolor}} 25 &{\cellcolor{shadecolor}} 37 &{\cellcolor{shadecolor}} -5.3834 &{\cellcolor{shadecolor}} 0.6545 &{\cellcolor{shadecolor}} -5.3834 &{\cellcolor{shadecolor}} 0.6545 &{\cellcolor{shadecolor}} -5.3834 &{\cellcolor{shadecolor}} 0.6545\tabularnewline
\hline 
26 & 27 & 2.573 & 0.6821 & 3.3466 & 0.7074 & 11.4067 & 2.0141\tabularnewline
\hline 
26 & 28 & -1.4082 & -0.2121 & -1.2867 & -0.2151 & -1.4543 & -0.7444\tabularnewline
\hline 
26 & 29 & -1.9019 & -0.2496 & -1.8111 & -0.2566 & -1.7398 & -0.9564\tabularnewline
\hline 
28 & 29 & -3.4761 & 0.2876 & -3.4761 & 0.2875 & -2.6488 & -1.0236\tabularnewline
\hline 
{\cellcolor{shadecolor}} 29 &{\cellcolor{shadecolor}} 38 &{\cellcolor{shadecolor}} -8.2477 &{\cellcolor{shadecolor}} 0.8033 &{\cellcolor{shadecolor}} -8.2477 &{\cellcolor{shadecolor}} 0.8033 &{\cellcolor{shadecolor}} -8.2477 &{\cellcolor{shadecolor}} 0.8033\tabularnewline
\hline 
\end{tabular}

\end{table*}

\section{Conclusion}
\label{sec:conclusion}
\mrn{This study delves into the impacts of False Data Injection (FDI) attacks on power systems, with a particular focus on the often-overlooked AC FDI attacks. By leveraging the PowerWorld simulator, we have presented a unified approach for assessing the impacts of optimal and arbitrary AC FDI attacks on power transmission lines' flows.
Our findings indicate that there is a critical trade-off between the residuals of state variables and the corresponding impacts of the AC FDI attacks, whether arbitrary or optimal. Specifically, optimal attacks lead to smaller alterations in the measurements, i.e., smaller attack values, and their residuals compared to arbitrary attacks. However, having small residuals comes at the cost of less severe impacts on the power system compared to arbitrary attacks.
Through the implementation and analysis of the IEEE 39-bus test systems, we have demonstrated the importance of considering this trade-off in the design of AC FDI attacks. The results underscore the need for a balanced approach in optimizing attack strategies to minimize detection while maximizing impact. This paper contributes to the broader understanding of AC FDI attacks and provides a foundation for future work aimed at enhancing the security and resilience of modern power systems against such sophisticated threats.}

\bibliographystyle{IEEEtran}
\IEEEtriggeratref{40}
\bibliography{ref}
\end{document}